\documentclass[12pt]{article}

\usepackage{amsmath,amssymb,enumerate,bbm}
\newtheorem{theorem}{Theorem}[section]

\newtheorem{lemma}[theorem]{Lemma}

\newtheorem{conjecture}[theorem]{Conjecture}
\numberwithin{equation}{section}
\newenvironment{proof}{\noindent\textbf{Proof\ }}{\hspace*{\fill}$\Box$\medskip}

 \makeatletter
 
 \newcommand{\Rmnum}[1]{\expandafter\@slowromancap\romannumeral #1@}
 \makeatother
\begin{document}

\title{On the spectrum of unitary finite-Euclidean graphs}\author{Si Li and Le Anh Vinh\footnote{Supported by the Ewrin Schr\"odinger International Institute for Mathematical Physics.}\\Mathematics Department\\Harvard University\\Cambridge MA 02138, US\\sili,vinh@math.harvard.edu}\maketitle

\begin{center}Mathematics Subject Classifications: 05C35, 05C38, 05C55, 05C25.\end{center}

\begin{abstract}
  We consider unitary graphs attached to $\mathbbm{Z}_n^d$ using
  an analogue of the Euclidean distance. These graphs are shown to be
  integral when $n$ is odd or the dimension $d$ is even. 
\end{abstract}

\section{Introduction}

Let $\Gamma$ be an additive group. For $S \subseteq \Gamma$, $0 \not \in S$ and
$S^{- 1} =\{- s : s \in S\}= S$, the Cayley graph $G = C (\Gamma, S)$ is the
undirected graph having vertex set $V (G) = \Gamma$ and edge set $E (G) =\{(a,
b) : a - b \in S\}$. The Cayley graph $G = C (\Gamma, S)$ is regular of degree
$|S|$. For a positive integer $n > 1$, the unitary Cayley graph $X_n = C (\mathbbm{Z}_n,
\mathbbm{Z}_n^{*})$ is defined by the additive group of the ring $\mathbbm{Z}_n$ of integers modulo
$n$ and the multiplicative group $\mathbbm{Z}_n^{*}$ of its units. So $X_n$ has vertex set $V
(X_n) = \mathbbm{Z}_n =\{0, 1, \ldots, n - 1\}$ and edge set \[E (X_n) =\{(a, b) : a, b
\in \mathbbm{Z}_n, \gcd (a - b, n) = 1\}.\] The graph $X_n$ is regular of degree $|\mathbbm{Z}_n^{*} |
= \phi (n)$, where $\phi (n)$ denotes the Euler function. Unitary Cayley
graphs are highly symmetric and have some remarkable properties connecting
graph theory and number theory. More information about the unitary Cayley
graphs can be found in Berrizbeitia and Giudici \cite{berrizbeitia}, Dejter and Giudici \cite{dejter},
Fuchs \cite{fuchs 1, fuchs 2}, and Klotz and Sander \cite{klotz}.

In this paper, we will study the higher dimensional unitary Cayley graphs, which
means we define unitary graphs over $\mathbbm{Z}_n^d$ for any integer $n$ and $d
\geqslant 1$. One possible generalization is to define the Cayley graph
$X_n^{(d)}$ with vertex set $V (X_n^{(d)}) = \mathbbm{Z}_n^d$ and edge
set
\begin{equation*}
  E (X_n^{(d)}) = \left\{ (a, b) : d (a, b) = \sum_{i = 1}^d (a_i - b_i) \in
  \mathbbm{Z}_n^{*} \right\} .
\end{equation*}
This generalization is however not very interesting as most of the known
results can be derived transparently. Thus, we consider unitary graphs
attached to $\mathbbm{Z}_n^d$ in a different way, using an analogue of the Euclidean
distance. Precisely, we define the unitary finite-Euclidean graphs $T_n^{(d)}$ with vertex
set $V (T_n^{(d)}) = \mathbbm{Z}_n^d$ and edge set
\begin{equation}
  E (T_n^{(d)}) = \left\{ (a, b) : d (a, b) = \sum_{i = 1}^d (a_i - b_i)^2 \in
  \mathbbm{Z}_n^{*} \right\} .
\end{equation}
Note that the finite-Euclidean graph $E_R^{(d)}(r)$ over a finite ring $R$ is the Cayley graph with vertex set $V(E_R^{(d)}(r)) = R^d$ and the edge set
\begin{equation*}
  E (E_R^{(d)}(r)) = \left\{ (a, b) : d (a, b) = \sum_{i = 1}^d (a_i - b_i)^2 = r \in R \right\} .
\end{equation*}

In \cite{medrano 1}, Medrano, Myers, Stark and Terras studied the spectrum of the finite-Euclidean graphs over finite fields and showed that these graphs are asymptotically Ramanujan graphs. In \cite{medrano 2}, these authors studied the same problem for the finite-Euclidean graphs over rings $\mathbbm{Z}_q$ for $q$ odd prime power. They showed that, over rings, except for the smallest case, the graphs (with unit distance parameter) are not (asymptotically) Ramanujan. 

In \cite{bannai}, Bannai, Shimabukuro and Tanaka showed that the finite-Euclidean graphs over finite fields are always asymptotically Ramanujan for a more general setting (i.e. they replace the Euclidean distance by nondegenerated quadratic forms). The second author recently applied these results to several interesting combinatorial problems, for example to tough Ramsey graphs (with P. Dung) \cite{dung}, Szemeredi-Trotter type theorem and sum-product estimate \cite{vinh 2} and the Erd\"os distance problem \cite{vinh 1}.

The main purpose of this paper is to study the spectrum of unitary finite-Euclidean graphs. We will show that the spectrum of unitary finite-Euclidean graphs consists entirely of integers when $n$ is odd or the dimension $d$ is even. This property seems to be amazingly widespread among Cayley graphs on abelian groups. One of the first papers in this direction is due to L. Lov\'asz \cite{lovasz 1}, who proved that all Cayley graphs, (cube-like)graphs, on $\mathbbm{Z}_2^d$ are integral where $\mathbbm{Z}_n$ is the ring of integers modulo n. See also W. So \cite{so} for related results.

The rest of this paper is organized as follows. In Section 2 we summarize preliminary facts that will be used throughout the paper. In Section 3 we prove the main theorem of this paper stating that unitary finite-Euclidean graphs are integral when $n$ is odd or the dimension $d$ is even. We conjecture that the result is also true for the odd dimensional case.

\section{Preliminaries}

We define the following conventions which will be used throughout the paper. 

\begin{itemize}
 \item The notation $\delta(\mathcal P)$
where $\mathcal P$ is some property, means that it's 1 if $\mathcal P$ is satisfied, 0 otherwise. 
 \item For any $k \in \mathbbm{Z}_n^{*}$ we define $I_n(k)$ to be the unique element $I_n(k) \in \mathbbm{Z}_n^{*}$ satisfying
\[k I_n(k)\equiv 1 \mod n.\]
 \item For any integer $k$ and $x=(x_1,\ldots,x_d) \in \mathbbm{Z}_n^{d}$, $k \mid x$ if and only if $k \mid x_i$ for all $1 \leq i \leq d$. 
 \item The exponential:
\[ e_n (x) = \exp \{2 \pi i x / n\}. \]
 \item We write $n = p_1^{r_1} \ldots p_t^{r_t}$ for the prime decomposition of
the positive integer $n$. For any nonempty subset $I \subseteq \{1, \ldots, t\}$, set $p_I =
\Pi_{i \in I} p_i$ and $n_I = n / p_I$. We define the integral square root of a positive integer $n$ to be
the largest integer $d$ such that for any $x$ satisfying $n \mid x^2$ then $d \mid x$. It is easy to shown that
if $n=p_1^{r_1}\ldots p_t^{r_t}$ then $d = p_1^{\lfloor \frac{r_1+1}{2} \rfloor} \ldots p_t^{\lfloor \frac{r_t+1}{2} \rfloor}$. We define $n_I^{'}$ to be the integral square root of $n_I$ and $p_I^{'} = n / n_I^{'}$. It is clear that $p_I \mid p_I^{'}$ and $n | n_Ip_I^{'}$ for any nonempty subset $I$ of $\{1,\ldots,t\}$.
 
\end{itemize}

In the following discussion, we will need to consider Gauss, Ramanujan and Gauss character sums over rings. 
For $c \in \mathbbm{Z}_n^{*}$ define the Gauss sum
\begin{equation}
  G_n (c) = \sum_{k \in \mathbbm{Z}_n} e_n (c k^2) .
\end{equation}

If $n$ odd, then
\begin{equation}
G_n(1) = \varepsilon_n \sqrt{n},
\end{equation}
where
\begin{equation*}
\varepsilon_n = \left\{ \begin{array}{cc}  
		1 & \mbox{if}\  n\equiv 1 (\mod 4)\\
		i & \mbox{if}\ n\equiv 3 (\mod 4).
	\end{array} \right.
\end{equation*}

If $(c,n)=1$ and $n$ odd, then 
\begin{equation}
G_n(c) = \left( \frac{c}{n} \right) G_n(1),
\end{equation}
here the Jacobi symbol $\left( \frac{c}{n} \right)$ is defined as follows:
\begin{itemize}
\item If $n = p$ is prime then 
\begin{equation*}
\left(\frac{c}{p}\right) = \left\{
        \begin{array}{ll}
        1 & p \nmid c, \ \ c \ \mbox{is square mod} \ p\\
        -1 & p \nmid c, \ \ c \ \mbox{is nonsquare mod} \ p\\
        0 & p \mid c.
    \end{array}  \right.
\end{equation*}
\item If $n = p_1^{r_1}\ldots p_t^{r_t}$ is the prime decomposition of $n$ then 
\begin{equation*}
\left(\frac{c}{n}\right) = \left(\frac{c}{p_1}\right)^{r_1}\ldots \left(\frac{c}{p_t}\right)^{r_t}.
\end{equation*}
\end{itemize}

If $(c,n)=1$ and $n \equiv 2 \mod 4$, then
\begin{equation} \label{4n+2}
G_n(c) = 0.
\end{equation}

If $(c,n)=1$ and $n \equiv 0 \mod 4$, then
\begin{equation} \label{4n}
G_n(c) = \left(\frac{n}{c}\right)(1+i^c)\sqrt{n}.
\end{equation}
For more information about the Gauss sums, we refer the reader to Section 1.5 in \cite{berndt}.

For $c \in \mathbbm{Z}_n$ define the Ramanujan sum
\begin{equation}
  r (c,n) = \sum_{(k, n) = 1} e_n (k c) .
\end{equation}
For positive integers $n, r$, Ramanujan sums have only integral values. Precisely, we have (cf. Corollary 2.4 of \cite{mccarthy})
\begin{equation}
r(c,n) = \mu(t_c) \frac{\phi(n)}{\phi(t_c)},
\end{equation}
where $t_c = \frac{n}{\gcd (c,n)}$. Here $\mu$ denotes the M\"obius function. 

For $c \in \mathbbm{Z}_n$ define the Gauss character sum for the quadratic character
\begin{equation}
  G_n (\chi, c) = \sum_{k \in \mathbbm{Z}_n} \left( \frac{k}{n} \right) e_n (c k).
\end{equation}
Then it can be shown that
\begin{equation} \label{gauss character} 
  G_n(\chi, c) = \left( \frac{c}{n} \right) G_n(\chi, 1).
\end{equation}  
If $n$ is square-free and $(c,n)=1$ then
\begin{equation} \label{gauss 1} 
G_n(c) = G_n(\chi,c).
\end{equation}

In general we have the following relation for the Gauss character sum and the Gauss sum.

\begin{lemma} \label{lemma 2.1}
  For any odd positive integer $n$ and $c \in \mathbbm{Z}_n$ we have $n \mid G_n(1) G_n (\chi, c)$. 
\end{lemma}

\begin{proof}
  From (\ref{gauss character}), it is suffice to prove for the case $c = -b^2$ for some $b \in \mathbbm{Z}_n^{*}$. We have
  \begin{eqnarray*}
  \sum_{k \in \mathbbm{Z}_n, x \in \mathbbm{Z}_n} e_n (b x + k x^2) & = & \sum_{x \in \mathbbm{Z}_n} e_n (b x) \sum_{k \in \mathbbm{Z}_n} e_n (x^2 k) \\
  & = &  \sum_{x \in \mathbbm{Z}_n} e_n(b x) \delta(n|x^2) n\\
  & = & nA, 
  \end{eqnarray*}   
  for some integer $A$. 
  Thus, by the Inclusion-Exclusion principle, we have
  \begin{eqnarray*}
  nA & = &  \sum_{k \in \mathbbm{Z}_n, x \in \mathbbm{Z}_n} e_n (b x + k x^2)\\
    & = & \sum_{\substack{(k, n) = 1\\x\in \mathbbm{Z_n}}} e_n(b x + k x^2) - \sum_{I \subseteq \{1, \ldots, t\}}(- 1)^{|I|}\sum_{\substack{s \in \mathbbm{Z}_{n_I}\\ x \in \mathbbm{Z}_n}} e_n (b x + p_I s x^2)\\
    & = & \sum_{\substack{(k,n) =1\\x\in \mathbbm{Z_n}}}e_n(k(x+I_n(2k)b)^2-I_n(4k)b^2) - \sum_{I \subseteq \{1, \ldots, t\}} (- 1)^{|I|} n_I \sum_{x \in \mathbbm{Z}_n} \delta (n_I |x^2_{}) e_n (bx)\\
    & = & \sum_{\substack{(k,n) =1\\x\in \mathbbm{Z_n}}}e_n(kx^2-I_n(4k)b^2) -  \sum_{I \subseteq \{1, \ldots, t\}} (- 1)^{|I|} n_I \sum_{s \in \mathbbm{Z}_{p_I'}} e_{p_I'} (b s)\\
    & = & \sum_{(k,n)=1} G_n(k)e_n(-I_n(4k)b^2) - \sum_{I \subseteq \{1, \ldots, t\}} (- 1)^{|I|} n_I \delta (p_I' |b) p_I'\\
    & = & \sum_{(k,n)=1} G_n(1) \left( \frac{k}{n} \right) e_n(-I_n(4k)b^2) - n B\\
    & = & G_n(1) \sum_{(k,n)=1}  \left( \frac{I_n(4k)}{n} \right) e_n(-I_n(4k)b^2) - n B \\
    & = & G_n(1) G_n (\chi, c) - n B,
  \end{eqnarray*}
  for some integer $B$ as $n|n_I p_I'$. The lemma follows.  
\end{proof}

\section{Unitary finite-Euclidean graphs}

Let 
\begin{equation}
S_d(n) := \{ x = (x_1,\ldots,x_d) \in \mathbbm{Z}_n^{d} | d(x,0) = \sum_{1}^{d}x_i^2 \in \mathbbm{Z}_n^{*} \}.
\end{equation}
Then the unitary finite-Euclidean graph $T_n^{(d)}$ is the Cayley graph $C(\mathbbm{Z}_n^d,S_d(n))$. Recall that the eigenvalues of Cayley graphs of abelian groups can be computed easily in terms of the characters of the group. This result, described in, e.g., \cite{lovazs}, implies that the eigenvalues of the unitary finite-Euclidean graph $T_n^{(d)}$ are all the numbers
\begin{equation}
\lambda_b = \sum_{x\in S_d(n)} e_n (^tb.x),
\end{equation}
where $b \in \mathbbm{Z}_n^{d}$. We will show that $\lambda_b$ is integer for any $b \in \mathbbm{Z}_n^{d}$ if $n$ is odd or $d$ is even.

For positive integers $n_1, n_2, n$ and $b \in \mathbbm{Z}_{n_1 n_2}^d$, we define
\begin{eqnarray*}
  f_{n_1, n_2} (b) & = & \sum_{k \in \mathbbm{Z}_{n_2}, x \in \mathbbm{Z}_{n_1 n_2}^d} e_{n_1 n_2} (^t
  b.x + n_1 k^t x.x) \\
  f_n (b) & = & \sum_{k \in \mathbbm{Z}_n, x \in \mathbbm{Z}_n^d} e_n (^t b.x +^{} k^t x.x)\\
  g_{n} (b) & = & \sum_{(k, n) = 1, x \in \mathbbm{Z}_n^d} e_n (^t b.x + k^t x.x) .
\end{eqnarray*}
  
  We first need some lemmas.
  
\begin{lemma} \label{lemma 3.1}
  For any $n_1, n_2$ positive integers and $b \in \mathbbm{Z}_{n_1 n_2}^d$ we have
  \begin{equation}
    f_{n_1, n_2} (b) = \delta (n_1|b) n_1^d f_{n_2} (c),
  \end{equation}
  where $c \in \mathbbm{Z}_{n_2}^d$ satisfies $b = n_1 c$ given that $\delta (n_1|b)$.
\end{lemma}

\begin{proof}
  We can write $x = n_2 x_1 + x_2$ uniquely where $x_1 \in \mathbbm{Z}_{n_1}^d$ and $x_2 \in
  \mathbbm{Z}_{n_2}^d$. Then
  \begin{eqnarray*}
    f_{n_1,n_2}(b) & = & \sum_{k \in \mathbbm{Z}_{n_2}, x \in \mathbbm{Z}_{n_1 n_2}^d} e_{n_1 n_2} (^t b.x + n_1 k ^tx.x) \\
    & = & \sum_{k \in \mathbbm{Z}_{n_2}, x_1 \in \mathbbm{Z}_{n_1}^d, x_2 \in \mathbbm{Z}_{n_2}^d} e_{n_1} (^t b.x_1) e_{n_1 n_2} (^t b.x_2 + n_1 k^t x_2 .x_2)\\
    & = & \delta (n_1|b) n_1^d \sum_{k \in \mathbbm{Z}_{n_2}, x_2 \in \mathbbm{Z}_{n_2}^d} e_{n_1n_2} (^t b.x_2 + n_1 k ^t x_2 .x_2)\\
    & = & \delta (n_1|b) n_1^d f_{n_2} (c) .
  \end{eqnarray*}
  This completes the proof of the lemma.
\end{proof}

\begin{lemma} \label{lemma 3.2}
  Suppose that $n = p_1^{r_1} \ldots p_t^{r_t}$ is the prime decomposition of a positive integer
  $n$. Then
  \begin{equation}
    f_n (b) = g_n (b) - \sum_{I \subseteq \{1, \ldots, t\}} (- 1)^{|I|}
    f_{p_I, n_I} (b) .
  \end{equation}
\end{lemma}

\begin{proof}
  By the Inclusion-Exclusion principle, we have
  \begin{eqnarray*}
    f_n(b) & = & \sum_{k \in \mathbbm{Z}_n, x \in \mathbbm{Z}_n^d} e_n (^t b.x + k^t x.x) \\
    & = & \sum_{(k, n) =
    1, x \in \mathbbm{Z}_n^d} e_n (^t b.x + k^t x.x) - \sum_{I \subseteq \{1, \ldots,
    t\}} (- 1)^{|I|} \sum_{s \in \mathbbm{Z}_{n_I}, x \in \mathbbm{Z}_n^d} e_n (^t b.x + p_I s^t x.x)\\
    & = & g_n (b) - \sum_{I \subseteq \{1, \ldots, t\}} (- 1)^{|I|} f_{p_I, n_I} (b),
  \end{eqnarray*}
  completing the lemma.
\end{proof}

\begin{lemma} \label{lemma 3.3}
  Suppose that $n = p_1^{r_1} \ldots p_t^{r_t}$ is the prime decomposition of
  a positive integer $n$. For any $b \in \mathbbm{Z}_n^d$ then
  \begin{equation}
    \lambda_b = \sum_{x \in \mathbbm{Z}_n^d : \gcd (^t x.x, n) = 1} e_n (^t b.x) =
    \delta (n|b) n^d + \sum_{I \subseteq \{1, \ldots, t\}} \frac{(-1)^{|I|}}{p_I} \delta(n_I|b)n_I^d f_{p_I}(b).
  \end{equation}
\end{lemma}

\begin{proof}
  By the Inclusion-Exclusion principle, we have
  \begin{eqnarray*}
    \lambda_b & = & \sum_{x \in \mathbbm{Z}_n^d} e_n (^t b.x) + \sum_{I \subseteq \{1,
    \ldots, t\}} (- 1)^{|I|} \sum_{^t x.x \in p_I \mathbbm{Z}_{n_I}} e_n (^t b.x)\\
    & = & \delta (n|b) n^d + \sum_{I \subseteq \{1, \ldots, t\}} \left\{
    \frac{(- 1)^{|I|}}{n} \sum_{k \in \mathbbm{Z}_{n_I}, h \in \mathbbm{Z}_n, x \in \mathbbm{Z}_n^d} e_n (^t b.x + h (^t x.x - p_I k)) \right\}\\
    & = & \delta (n|b) n^d + \sum_{I \subseteq \{1, \ldots, t\}} \left\{
    \frac{(- 1)^{|I|}}{n} \sum_{x \in \mathbbm{Z}_n^d, h \in \mathbbm{Z}_n} \left( e_n (^t b.x +
    h ^tx.x) \sum_{k \in \mathbbm{Z}_{n_I}} e_{n_I} (- h k) \right) \right\}\\
    & = & \delta (n|b) n^d + \sum_{I \subseteq \{1, \ldots, t\}} \left\{
    \frac{(- 1)^{|I|}}{p_I} \sum_{x \in \mathbbm{Z}_n^d, h \in n_I \mathbbm{Z}_{p_I}} e_n (^t b.x +
    h^t x.x) \right\}\\
    & = & \delta (n|b) n^d + \sum_{I \subseteq \{1, \ldots, t\}} \frac{(-
    1)^{|I|}}{p_I} f_{n_I, p_I} (b)\\
    & = & \delta (n|b) n^d + \sum_{I \subseteq \{1, \ldots, t\}} \frac{(-
    1)^{|I|}}{p_I} \delta(n_I|b)n_I^d f_{p_I}(b),
  \end{eqnarray*}
  where the last line follows from Lemma \ref{lemma 3.1}.
  This concludes the proof of the lemma. 
\end{proof}

If $n$ is odd, we can write $g_n(b)$ in terms of the Gauss and Ramanujan sums.

\begin{lemma} \label{lemma 3.4}
  a) For any odd positive integer $n$ we have
  \begin{equation}
  g_n(b) = \left\{ \begin{array}{ll} 
    (G_n(1))^d G_n(\chi,- ^t b.b) & d\ \mbox{odd}\\
    (G_n(1))^d r(-^t b.b, n) & d\ \mbox{even}.
    \end{array} \right.
  \end{equation}
  b) Furthermore, $n \mid g_n(b)$ for all odd integers $n$.
\end{lemma}

\begin{proof}
  a) We have
  \begin{eqnarray*}
    \sum_{\substack{(k, n) = 1\\ x \in \mathbbm{Z}_n^d}} e_n (^t b.x + k^t x.x) & = & \sum_{(k, n) = 1, x \in \mathbbm{Z}_n^d} e_n (k^t(x + I_n (2 k)b).(x + I_n(2k)b)- I_n(4k)^t b.b)\\
    & = & \sum_{(k, n) = 1, x \in \mathbbm{Z}_n^d} e_n (k^t x.x - I_n (4 k)^t b.b)\\
    & = & \sum_{(k, n) = 1} G_n ^d(k) e_n (- I_n (4 k)^t b.b)\\
    & = & (G_n(1))^d \sum_{(k, n) = 1} \left( \frac{k}{n}
    \right)^d e_n (- I_n (4 k)^t b.b)\\
    & = & (G_n(1))^d \sum_{(k, n) = 1} \left( \frac{I_n (4
    k)}{n} \right)^d e_n (- I_n (4 k)^t b.b)\\
    & = & (G_n(1))^d \sum_{(k, n) = 1} \left( \frac{k}{n}
    \right)^d e_n (-^t b.b k)\\
    & = & \left\{ 
    \begin{array}{ll} 
    (G_n(1))^dG_n(\chi,- ^t b.b) & d\ \mbox{odd}\\
    (G_n(1))^dr(-^t b.b, n) & d\ \mbox{even}
    \end{array} \right.
  \end{eqnarray*}
  This completes the first part of the lemma.
  
  b) We have
  \begin{eqnarray*}
    g_n(b) & = & \left\{ 
    \begin{array}{ll} 
    (G_n(1))^dG_n(\chi,- ^t b.b) & d\ \mbox{odd}\\
    (G_n(1))^dr(-^t b.b, n) & d\ \mbox{even}
    \end{array} \right.\\
    & = & \left\{ 
    \begin{array}{ll} 
    \varepsilon_n^{d-1}n^{(d-1)/2}G_n(1)G_n(\chi,- ^t b.b) & d\ \mbox{odd}\\
    \varepsilon_n^{d}n^{d/2}r(-^t b.b, n) & d\ \mbox{even}
    \end{array} \right. \\
  \end{eqnarray*}
  Part b) follows immediately from Lemma \ref{lemma 2.1} and the fact that the Ramanujan sums have only integral values. 
\end{proof}

We also can write $g_n(b)$ in form of the Gauss sums when $n$ is even.

\begin{lemma} \label{lemma 3.7}
a) Suppose that $n$ is a squarefree even integer and $b = \{b_i\}_1^n \in \mathbbm{Z}_n^d$, then
  \begin{equation*}
  g_n(b) = \left\{ \begin{array}{ll}
  0 & \mbox{there exists} \ i \ \mbox{such that} \ 2|b_i\\
  \sum_{(k, n) = 1} e_{4 n} (- I_{4 n} (k)^t b.b) (G_{4 n} (k)
    / 2)^d & 2 \nmid b_i \ \mbox{for all} \ i.
  \end{array}
  \right.
  \end{equation*}

b) Suppose that $2 \mid d$. For any squarefree even integer $n$ and any vector $b \in Z_n^d$ we have $n \mid g_n(b)$.
\end{lemma}

\begin{proof}
a) We have
  \begin{eqnarray}
    g_n (b) & = & \sum_{(k, n) = 1} \sum_{x \in Z_n^d} e_n (^t b.x + k^t x.x)
    \nonumber\\
    & = & \sum_{(k, n) = 1} \prod_{i = 1}^d \left( \sum_{x \in Z_n} e_n (b_i
    x + k x^2) \right) . \nonumber
  \end{eqnarray}
  Suppose that there exists $i$ such that $2| b_i$. If we write $b_i = 2 c$ then
  \begin{eqnarray}
    \sum_{x \in Z_n} e_n (b_i x + k x^2) & = & \sum_{x \in Z_n} e_n (k (x + c
    I_n (k))^2 - I_n (k) c^2) \nonumber\\
    & = & e_n (- I_n (k) c^2) G_n (k) \nonumber\\
    & = & 0, \nonumber
  \end{eqnarray}
  where the last line follows from (\ref{4n+2}). Thus, $g_n (b) = 0$ if there exists $i$ such that $2| b_i$. Now, suppose
  that $2 \nmid b_i$ for all $1 \leqslant i \leqslant d$. We have
  \begin{eqnarray}
    g_n (b) & = & \sum_{(k, n) = 1, x \in Z_n^d} e_{4 n} (4^t b.x + 4 k^t x.x)
    \nonumber\\
    & = & \sum_{(k, n) = 1, x \in Z_n^d} e_{4 n} (k^t (2 x + I_{4 n} (k) b) .
    (2 x + I_{4 n} (k) b) - I_{4 n} (k)^t b.b) \nonumber\\
    & = & \sum_{(k, n) = 1} e_{4 n} (- I_{4 n} (k)^t b.b) \left\{ \prod_{i =
    1}^d \left( \sum_{x \in Z_n} e_{4 n} (k (2 x + I_{4 n} (k) b_i)^2 \right)
    \right\} . \nonumber
  \end{eqnarray}
  Let $a_i = I_{4 n} (k) b_i$, then $a_i$ is odd. Set
  \begin{equation}
    S_i = \sum_{x \in Z_n} e_n (k (2 x + a_i)^2) .
  \end{equation}
  Substitute $x = x + n$ into $S_i$, we have
  \begin{eqnarray}
    S_i & = & \frac{1}{2} \sum_{x \in Z_{2 n}} e_{4 n} (k (2 x + a_i)^2)
    \nonumber\\
    & = & \frac{1}{2} \left( \sum_{x \in Z_{4 n}} e_{4 n} (k x^2) - \sum_{x
    \in 2 Z_{2 n}} e_{4 n} (k x^2) \right) \nonumber\\
    & = & \frac{1}{2} \left( G_{4 n} (k) - \sum_{x \in Z_{2 n}} e_n (k x^2)
    \right) \nonumber\\
    & = & \frac{1}{2} (G_{4 n} (k) - 2 G_n (k)) = \frac{1}{2} G_{4 n} (k),
    \nonumber
  \end{eqnarray}
  because $G_n (k) = 0$ for $n \equiv 2$ mod $4$ and $(k, n) = 1$. Therefore, we
  have
  \begin{equation*}
    g_n (b) =  \sum_{(k, n) = 1} e_{4 n} (- I_{4 n} (k)^t b.b) (G_{4 n} (k)
    / 2)^d .
  \end{equation*}
  This concludes the proof of part a). 

b)  If there exists $i$ such that $2 | b_i$ then from part a), $g_n(b) = 0$. Suppose that $2 \nmid b_i$ for all $i$.
  Substitute (\ref{4n}) into part a). Note that $k$ odd so $(1 + i^k)^2 = 2 i^k$ and $(1 + i^k)^4 = - 4$. We consider
  two cases.
  
  Case 1: Suppose that $d = 4 d_1$. Since $2 \nmid b_i$ for all $i$, we have
  $4 \mid^t b.b$.
  \begin{eqnarray}
    g_n (b) & = & \sum_{(k, n) = 1} e_{4 n} (- I_{4 n} (k)^t b.b) (1 + i^k)^{4
    d_1} n^{2 d_1} \nonumber\\
    & = & (- 4)^{d_1} n^{2 d_1} \sum_{(k, n) = 1} e_n \left( \frac{-^t
    b.b}{4} I_n (k) \right) \nonumber\\
    & = & (- 4)^{d_1} n^{2 d_1} r (-^t b.b / 4, n) . \label{case 1}
  \end{eqnarray}
  
  Case 2: Suppose that $d = 4 d_1 + 2$. Since $2 \nmid b_i$ for all $i$, we
  have $-^t b.b = 2 c$ for some odd $c$. Let $n = 2 m$ for some $m$ odd. We
  have
  \begin{eqnarray}
    g_n (b) & = & \sum_{(k, n) = 1} e_{2 n} (I_{2 n} (k) c) (1 + i^k)^{4 d_1 +
    2} n^{2 d_1 + 1} \nonumber\\
    & = & (- 1)^{d_1} (2 n)^{2 d_1 + 1} \sum_{(k, n) = 1} e_{2 n} (I_{2 n}
    (k) c) i^k \nonumber\\
    & = & (- 1)^{d_1} (2 n)^{2 d_1 + 1} \sum_{(k, 2 m) = 1} e_{4 m} ((I_{4 m}
    (k) c) i^{I_{4 m} (k)} \nonumber\\
    & = & (- 1)^{d_1} (2 n)^{2 d_1 + 1} \sum_{(k, 2 m) = 1} e_{4 m} (I_{4 m}
    (k) (c + m)) \nonumber\\
    & = & (- 1)^{d_1} (2 n)^{2 d_1 + 1} \sum_{(k, 2 m) = 1} e_{2 m} (I_{2 m}
    (k) (c + m) / 2) \nonumber\\
    & = & (- 1)^{d_1} (2 n)^{2 d_1 + 1} r ((c + m) / 2, n) \label{case 2}
  \end{eqnarray}
  since $I_{4 m} (k) \equiv k$ mod $4$ and $I_{4 m} (k) \equiv I_{2 m} (k)$
  mod $2 m$.
  
  Part b) follows immediately from (\ref{case 1}), (\ref{case 2}) and the fact that Ramanujan sums have only integral values.
\end{proof}

We are now ready to prove the main result of the paper. 

\begin{theorem} \label{main 1}
  Suppose that $n$ is an odd integer or $d$ is an even integer then all eigenvalues of the unitary finite-Euclidean graph $T_n^{(d)}$ are integers. 
\end{theorem}

\begin{proof}
  We consider two cases.
  
  Case 1: Suppose that $n$ is odd. From Lemmas \ref{lemma 3.1}, \ref{lemma 3.2} and \ref{lemma 3.4}, we can show by induction that $n \mid f_n (b)$ for any
  positive odd integer $n$ and for any vector $b \in \mathbbm{Z}_n^d$. Together with Lemma
  \ref{lemma 3.3}, we have $\lambda_b$'s are integer for all $b \in \mathbbm{Z}_n^d$. 
  
  Case 2: Suppose that $d$ is even.  From Lemmas \ref{lemma 3.2}, \ref{lemma 3.4} and \ref{lemma 3.7}, we can show by induction that $n \mid f_n(b)$ for any squarefree positive integer $n$ and for any vector $b \in \mathbbm{Z}_n^d$. Since $p_I$ is squarefree for all $I \subseteq \{1,\ldots,n\}$, from Lemma \ref{lemma 3.3}, we have $\lambda_b$'s are integer for all $b \in \mathbbm{Z}_n^d$. This completes the proof of the theorem.
\end{proof}

We conjecture that the same result also holds for odd dimensional cases. 

\begin{conjecture}
 For any positive integers $n$ and $d$, then all eigenvalues of the unitary finite-Euclidean graph $T_n^{(d)}$ are integers. 
\end{conjecture}

From Theorem \ref{main 1}, the remaining open case is: $n$ even and $d$ odd.


\begin{thebibliography}{99}

\bibitem{angle} J. Angle, B. Shook, A. Terras, C. Trimble and E. Velasquez, Graph spectra for finite upper half planes over rings, 
\textit{Linear Algebra Applications}, \textbf{226-228} (1995), 423--457.

\bibitem{bannai}
E. Bannai, O. Shimabukuro and H. Tanaka, Finite Euclidean graphs and Ramanujan graphs, \textit{Discrete Mathematics} (to appear).

\bibitem{berrizbeitia}
P. Berrizbeitia and R. E. Giudici, On cycles in the sequence of unitary Cayley graphs, \textit{Discrete Mathematics} \textbf{282} 1-3 (2004), 239--243.


\bibitem{berndt}
B. C. Berndt, R. J. Evans and K. S. Williams, \textit{Gauss and Jacobi Sums}, Canadian Mathematical Society Series of Monographs and Advanced Texts \textbf{21}, John Wiley \& Sons (1998).


\bibitem{dejter}
I. J. Dejter and R. E. Giudici, On unitary Cayley graphs, \textit{J. Combin. Math. Combin. Comput.}, \textbf{18} (1995), 121--124.

\bibitem{fuchs 1}
E. D. Fuchs, Longest induced cycles in circulant graphs, \textit{The Electronic Journal of Combinatorics}, \textbf{12} (2005), 1-12. 

\bibitem{fuchs 2}
E. D. Fuchs and J. Sinz, Longest induced cycles in Cayley graphs,\\eprint arXiv:math/0413308. 

\bibitem{klotz}
W. Klotz and T. Sander, Some properties of unitary Cayley graphs, \textit{The Electronic Journal of Combinatorics}, \textbf{14} (2007), R45. 

\bibitem{lovasz 1}
L. Lov\'asz, Spectra of graphs with transitive groups, \textit{Periodica Mathematica Hungarica}, \textbf{6} (1975), 191-195.

\bibitem{lovazs}
L. Lov\'asz, \textit{Combinatorial Problems and Exercises}, North Holland, Amsterdam, 1979, Problem 11.8.

\bibitem{mccarthy}
P. J. McCarthy, \textit{Introduction to arithmetical functions}, Universitext. Springer-Verlag, New York, 1994.

\bibitem
{medrano 1} A. Medrano, P. Myers, H. M. Stark and A. Terras, Finite analogues of Euclidean space,
\textit{Journal of Computational and Applied Mathematics}, \textbf{68} (1996), 221--238.

\bibitem
{medrano 2} A. Medrano, P. Myers, H. M. Stark and A. Terras, Finite Euclidean graphs over rings,
\textit{Proceedings of the American Mathematics Society}, \textbf{126}  (1998), 701--710.

\bibitem{so}
W. So, Integral circulant graphs, \textit{Discrete Mathematics}, \textbf{306} (2006), 153-158.

\bibitem
{dung}
L. A. Vinh and D. P. Dung, Explicit tough Ramsey graphs, 
\textit{Proceedings of the International Conference on Relations, Orders and Graphs: Interaction with Computer Science} 2008, Nouha Editions, 139-146. 

\bibitem{vinh 1}
L. A. Vinh, Explicit Ramsey graphs and Erd\"os distance problem over finite Euclidean and non-Euclidean spaces, \textit{The Electronic Journal of Combinatorics}, \textbf{15} (1), 2008, R5. 

\bibitem{vinh 2}
L. A. Vinh, Szemeredi-Trotter type theorem and sum-product estimate in finite fields, submitted.

\end{thebibliography}
\end{document}